\documentclass[10pt]{amsart}

\makeatletter
\@namedef{subjclassname@2020}{\textup{2020} Mathematics Subject Classification}
\makeatother

\usepackage{amsfonts, amsmath, amssymb, enumerate, lineno, mathrsfs}

\newtheorem{theorem}{Theorem}[section]
\newtheorem{lemma}[theorem]{Lemma}

\newtheorem{corollary}[theorem]{Corollary}
\theoremstyle{definition}
\newtheorem{definition}[theorem]{Definition}
\numberwithin{equation}{section}

\def\cf{\mathop{\rm cf}\nolimits}
\def\dom{\mathop{\rm dom}\nolimits}

\let\mathcal\mathscr

\begin{document}


\title{Singular precalibers for topological products}

\author[A. R\'{\i}os-Herrej\'on]{Alejandro R\'{\i}os-Herrej\'on}
\address{Departamento de Matem\'aticas, Facultad de Ciencias, Circuito ext. s/n, Ciudad Universitaria, C.P. 04510,  M\'exico, CDMX}
\email{chanchito@ciencias.unam.mx}
\urladdr{}
\thanks{The author was supported by CONACYT grant no. 814282}



\begin{abstract}
In this paper we prove that if $\kappa$ is a singular cardinal with uncountable cofinality, then every power of a given topological space with precaliber $\kappa$ has precaliber $\kappa$ as well. Furthermore, if $\{X_\alpha : \alpha<\lambda\}$ is a family of topological spaces with precaliber $\kappa$ and $\lambda<\cf(\kappa)$, then $\kappa$ is a precaliber for the topological product $\prod\{X_\alpha : \alpha<\lambda\}$.
\end{abstract}

\maketitle

\section{Introduction}

The concepts of ``caliber'' and ``precaliber'' for a topological space were established by \v{S}anin in his papers \cite{sanin1946a}--\cite{sanin1948} between 1946 and 1948. Since then, the mathematical production that has emerged from these notions has been extensive and enlightening (see, for instance, \cite{mcintyre2006}, \cite{shakhmatov1986} and \cite{tall1977}--\cite{wathao1990}).

For example, in the classical book \cite{comneg1982}, Comfort and Negrepontis study the natural generalizations of these concepts to pairs and triples of cardinal numbers. Furthermore, they explore some variations of the concept of caliber that appear naturally when working with these notions (compact-calibers, pseudo-compactness numbers, the Knaster property $(K)$, the Souslin number, etc.).

Additionally, there are close relationships between the concepts of precalibers and calibers with some cardinal functions, for instance, the cellularity, density and $\pi$-weight of a topological space. These facts provide tools that facilitate the study of certain objects when viewed from the unique perspective that comes from being familiar with the basic techniques and results of this area.

One of the first steps in developing this new machinery was trying to understand the behavior of calibers and precalibers with respect to topological products. A well-known result regarding this problem was proven by \v{S}anin himself in \cite{sanin1948} and states the following:

\begin{theorem}\label{thm_Sanin} Let $\kappa$ be an uncountable regular cardinal and $\{X_\alpha : \alpha<\lambda\}$ be a family of topological spaces. Then $\kappa$ is a precaliber (resp., caliber) for the topological product $\prod\{X_\alpha : \alpha<\lambda\}$ if and only if $\kappa$ is a precaliber (resp., caliber) for each of its factors.
\end{theorem}

In light of this result, a natural question arises: what happens when $\kappa$ is a singular cardinal? Over the years, there have been some partial advances on this issue. For example, in \cite{argtsa1982} Argyros and Tsarpalias showed that, under the Generalized Continuum Hypothesis, the following result is true:

\begin{theorem}\label{thm_Arg_Tsa} [\textsf{GCH}] If $\{X_\alpha : \alpha<\lambda\}$ is a family of topological spaces and $\kappa$ is a cardinal number with $\cf(\kappa)>\omega$, then the following statements are satisfied.
\begin{enumerate}
\item $\kappa$ is a precaliber for the topological product $\prod\{X_\alpha : \alpha<\lambda\}$ if and only if $\kappa$ is a precaliber for each of its factors.
\item If every $X_\alpha$ is a compact Hausdorff space, then $\kappa$ is a caliber for $\prod\{X_\alpha : \alpha<\lambda\}$ if and only if $\kappa$ is a caliber for each of its factors.
\end{enumerate}
\end{theorem}

A theorem of this kind, while nice to have at our disposal, is always accompanied by the question of whether or not it will be possible to prove a similar result within $\mathsf{ZFC}$. In this sense, a partial answer without using additional axioms was given by Shelah in \cite{shelah1977} with the next result:

\begin{theorem}\label{thm_Shelah} The following statements hold for any singular cardinal $\kappa$.

\begin{enumerate}
\item If $X$ is a topological space with caliber $\kappa$, $\cf(\kappa)>\omega$ and $\lambda$ is a cardinal number, then $X^{\lambda}$ also has caliber $\kappa$.
\item Whenever $\{X_\alpha : \alpha<\lambda\}$ is a family of topological spaces, each of them with caliber $\kappa$, and $\lambda < \cf(\kappa)$, then $\prod\{X_\alpha : \alpha<\lambda\}$ has caliber $\kappa$ as well.
\end{enumerate}
\end{theorem}

Now, the fact that Theorems~\ref{thm_Sanin} and \ref{thm_Arg_Tsa} have a version for precalibers and calibers but Theorem~\ref{thm_Shelah} only speaks about calibers leaves us with a little concern: could it be that the precaliber version of this statement is not true?

The purpose of this article is to show that, if we change the word ``caliber'' to ``precaliber'' in Theorem~\ref{thm_Shelah}, then the corresponding result is also correct. Furthermore, it is worth mentioning that, although the proof in this text follows mostly the original ideas of \cite{shelah1977}, certain key results are simplified and explained meticulously.

\section{Preliminaries}\label{secc_prelim}

All topological notions and all set-theoretic concepts whose definition is not included here should be understood as in
\cite{engelking1989} and \cite{kunen2013}, respectively.

The symbol $\omega$ will stand for both, the set of all non-negative integers and the first infinite cardinal. 

If $\kappa$ is a cardinal number, the {\it cofinality} of $\kappa$, denoted by $\cf(\kappa)$, will be the least ordinal for which there is a sequence of cardinal numbers $\{\kappa_\alpha : \alpha<\cf(\kappa)\}$ such that $\kappa=\sup\{\kappa_\alpha : \alpha<\cf(\kappa)\}$. Additionally, if $X$ is a set, we will denote the families $\{Y\subseteq X : |Y|<\omega\}$ and $\{Y\subseteq X : |Y|=\kappa\}$ by $[X]^{<\omega}$ and $[X]^{\kappa}$, respectively.

For a topological space $X$, we will use the symbol $\tau_X$ to refer to the family of open subsets of $X$. Similarly, $\tau_X^{+}$ will be used to denote the set $\tau_X \setminus \{\emptyset\}$.

A non-empty subset $\mathcal{U}$ of $\tau_X^+$ will be called a {\it centered family} if and only if every element of $[\mathcal{U}]^{<\omega}\setminus\{\emptyset\}$ has non-empty intersection. In particular, if $\{U_\gamma : \gamma\in I\}$ is a non-empty indexed subset of $\tau_X^+$, then $\{U_\gamma : \gamma\in I\}$ is a centered family if and only if for every $J\in [I]^{<\omega}\setminus\{\emptyset\}$, $\bigcap\{U_\gamma : \gamma\in J\}\neq\emptyset$.

Throughout the text and unless explicitly stated otherwise, all spaces and cardinal numbers considered will be infinite.

Let us begin with the following definition.

\begin{definition}\label{def} Let $X$ be a topological space and $\kappa$ be a cardinal number.

\begin{enumerate}
\item We will say that $X$ has {\it precaliber $\kappa$} (or rather, that $\kappa$ is a {\it precaliber for $X$}) if for any set $I$ of cardinality $\kappa$ and for any family $\{U_\gamma : \gamma\in I\} \subseteq \tau_{X}^{+}$, there exists $J\in [I]^{\kappa}$ such that $\{U_\gamma : \gamma\in  J\}$ is a centered family.
\item If $\lambda$ is a cardinal number that satisfies the relation $\lambda\leq \kappa$, we say that $X$ has {\it precaliber $(\kappa,\lambda)$} (or else, that $(\kappa ,\lambda)$ is a {\it precaliber} for $X$) if and only if for any set $I$ of cardinality $\kappa$ and for any family $\{U_\gamma : \gamma\in I \} \subseteq \tau_{X}^{+}$, there is $J\in [I]^{\lambda}$ such that $\{U_\gamma : \gamma\in J\}$ is a centered family.
\end{enumerate}

\end{definition}

Clearly, by establishing adequate bijective functions, it can be verified that a cardinal $\kappa$ (respectively, a pair of cardinals $(\kappa,\lambda)$) is a precaliber for a topological space $X$ if and only if for every family $\{U_\gamma : \gamma <\kappa\}$ contained in $\tau_X^+$ there is $J\in [\kappa]^{\kappa}$ (resp., $J\in [\kappa]^{\lambda}$) such that $\{U_\gamma : \gamma \in J\}$ is centered.

On the other hand, if $\mathcal{B}$ is a basis for $X$, then a routine argument shows that $\kappa$ (resp., $(\kappa,\lambda)$) is a precaliber for $X$ if and only if for every family $\{B_\gamma : \gamma <\kappa\}\subseteq \mathcal{B}\setminus\{\emptyset\}$ there is $J\in [\kappa]^{\kappa}$ (resp., $J\in [\kappa]^{\lambda}$) such that $\{B_\gamma : \gamma \in J\}$ is centered. This fact will be used in the proof of Theorem~\ref{central_thm}.

If $\kappa$ and $\lambda$ are cardinal numbers such that $\kappa<\lambda$, then the symbol $[\kappa,\lambda)$ will stand for the set of all ordinal numbers $\gamma$ that verify the inequalities $\kappa\leq \gamma<\lambda$.

In the remainder of this portion of the text we will present some auxiliary and known propositions that we will use in the proofs of the central theorems of section~\ref{secc_main}. For completeness purposes, we will add the arguments for some key results.

The next two results are straightforward. The first of these can be found in \cite[Corollary~2.25, p.~40]{comneg1982}, while the second one follows directly from Definition~\ref{def}.

\begin{lemma}\label{lemma_prec_cof} If $X$ is a topological space and $\kappa$ is a precaliber for $X$, then $\cf(\kappa)$ is also a precaliber for $X$.

\end{lemma}

\begin{lemma}\label{lemma_basic_pairs} The following statements are true for any topological space $X$.

\begin{enumerate}
\item $\kappa$ is a precaliber for $X$ if and only if $(\kappa,\kappa)$ is a precaliber for $X$.

\item Whenever $(\lambda,\mu)$ is a precaliber for $X$ and $\lambda' \geq \lambda \geq \mu \geq \mu'$, then $(\lambda',\mu')$ is also a precaliber for $X$.

\item If $\kappa > \mu \geq \cf(\mu) > \cf(\kappa)$ and, for every regular $\kappa>\lambda\geq \mu$, $(\lambda,\mu)$ is not a precaliber for $X$, then $(\kappa,\mu)$ is not a precaliber for $X$ either.
\end{enumerate}

\end{lemma}

\begin{definition}\label{def_P(kappa,X)} Let $\kappa$ be a singular cardinal and $X$ be a topological space. We will say that a collection $\{\kappa_\alpha : \alpha<\cf(\kappa)\}$ satisfies the $P(\kappa, X)$ condition if the following properties hold:
\begin{enumerate}
\item $\sup\{\kappa_\alpha : \alpha<\cf(\kappa)\}=\kappa$;
\item for each $\alpha<\cf(\kappa)$, $\kappa_\alpha$ is a regular and uncountable cardinal such that $\max\left\{\cf(\kappa),\sup\{\kappa_\beta : \beta<\alpha\}\right\}<\kappa_\alpha$; and
\item $(\kappa_{\alpha+1},\kappa_\alpha)$ is a precaliber for $X$, provided that $\alpha<\cf(\kappa)$.
\end{enumerate}

\end{definition}

\begin{lemma}\label{lemma_P(kappa,X)} If a singular cardinal $\kappa$ is a precaliber for a topological space $X$, then the following statements hold.

\begin{enumerate}
\item For any $\mu<\kappa$, there exists a regular cardinal $\mu< \lambda<\kappa$ such that $(\lambda,\mu)$ is a precaliber for $X$.

\item There is a collection $\{\kappa_\alpha : \alpha<\cf(\kappa)\}$ that satisfies the $P(\kappa, X)$ condition.

\end{enumerate}

\end{lemma}

\begin{proof} To argue (1) tet us fix $\mu<\kappa$. Note that the singularity of $\kappa$ ensures the relation $\mu^{+}<\kappa$. Now, if $\cf(\kappa)\geq \mu^{+}$, then by letting $\lambda:= \cf(\kappa)$, Lemma~\ref{lemma_prec_cof} guarantees that $X$ has precaliber $\lambda$ and, therefore, item (2) of Lemma~\ref{lemma_basic_pairs} tells us that $(\lambda,\mu^{+})$ is a precaliber for $X$. On the other hand, when $\mu^{+}>\cf(\kappa)$, we use a combination of items (2) and (3) of Lemma~\ref{lemma_basic_pairs} to find a regular cardinal $\kappa>\lambda \geq \mu^{+}$ such that $(\lambda,\mu^{+})$ is a precaliber for $X$. In any case, item (2) of Lemma~\ref{lemma_basic_pairs} guarantees that there is a regular cardinal $\mu< \lambda<\kappa$ such that $(\lambda,\mu)$ is a precaliber for $X$.

Now, with item (2) in mind, let $\{\lambda_\alpha : \alpha<\cf(\kappa)\}$ be a strictly increasing sequence of cardinal numbers such that $\sup\{\lambda_\alpha : \alpha<\cf(\kappa)\}=\kappa$.

\medskip

\noindent {\bf Claim.} There is a collection $\{\kappa_\alpha : \alpha<\cf(\kappa)\}$ that satisfies the following conditions for each $\alpha<\cf(\kappa)$:

\begin{enumerate}
\item $\max\left\{\cf(\kappa),\lambda_\alpha, \sup\{\kappa_\beta : \beta<\alpha\}\right\} < \kappa_\alpha<\kappa$; 
\item $\kappa_\alpha$ is regular and uncountable; and
\item $(\kappa_{\alpha+1},\kappa_\alpha)$ is a precaliber for $X$.
\end{enumerate}

\medskip

Indeed, by recursion suppose that for some $\alpha<\cf(\kappa)$ we have managed to construct $\{\kappa_\beta : \beta<\alpha\}$ with the desired properties. In these circumstances, if we let $\mu:= \max\left\{\cf(\kappa),\lambda_\alpha, \sup\{\kappa_\beta : \beta<\alpha\} \right\}$, then item (1) of the present Lemma gives us a regular cardinal $\mu < \kappa_\alpha < \kappa$ such that $\left(\kappa_{\alpha}, \mu\right)$ is a precaliber for $X$. In particular,  item (2) of Lemma~\ref{lemma_basic_pairs} implies that, if $\alpha=\gamma+1$, then $(\kappa_{\gamma+1},\kappa_\gamma) $ is a precaliber for $X$.

Finally, it is clear that the conditions (1) to (3) of our claim guarantee that $\{\kappa_\alpha : \alpha<\cf(\kappa)\}$ satisfies the $P(\kappa, X)$ condition.
\end{proof}

The following result can be proven with a routine argument and the next two corollaries follow directly from it.

\begin{lemma}\label{lemma_chain_prec} Let $X$ and $Y$ be a pair of topological spaces. If $\kappa$, $\lambda$ and $\mu$ are cardinals such that $\kappa\geq \lambda \geq \mu$, $(\kappa,\lambda)$ is a precaliber for $Y$ and $ (\lambda,\mu)$ is a precaliber for $X$, then $(\kappa,\mu)$ is a precaliber for $X\times Y$.

\end{lemma}

\begin{corollary}\label{cor_1} Let $n$ be a positive integer, $\{X_m : m< n\}$ be a family of topological spaces and $\{\kappa_m : m\leq n \}$ be an increasing sequence of cardinal numbers. If for each $m< n$ it is satisfied that $(\kappa_{m+1},\kappa_m)$ is a precaliber for $X_m$, then $(\kappa_{n} , \kappa_0)$ is a precaliber for $\prod\{X_m : m< n\}$.

\end{corollary}

\begin{corollary}\label{cor_2} Let $n$ be a positive integer, $\{X_m : m <n\}$ be a family of topological spaces and $\kappa$ be a cardinal number. If $\kappa$ is a precaliber for each of the $X_m$, then $\kappa$ is a precaliber for $\prod\{X_m : m< n\}$.
\end{corollary}

It only remains to mention a couple of set-theoretic results that will be the cornerstone of Theorem~\ref{central_thm}. The first result of this kind is a well-known consequence of the $\Delta$-system Lemma. A proof of this statement can be found in \cite[Lemma~\textsc{III}.2.6, p.~166]{kunen2013}.

\begin{theorem}\label{delta_system} If $\kappa$ is a cardinal number such that $\kappa=\cf(\kappa)>\omega$ and $\{A_\alpha : \alpha<\kappa\}$ is a family of finite sets, then there are sets $J\in [\kappa]^{\kappa}$ and $A$ such that $A_\alpha \cap A_\beta = A$, provided that $\alpha, \beta \in J$ are distinct.

\end{theorem}

In the following theorem we prove a more detailed version of Claim~2.4.(B) from \cite{shelah1977}. Apparently, based on the comments on page 19 of \cite{comneg1982}, this result has been proven in a couple of somewhat innacesible documents, namely, in the doctoral dissertation \cite{argyros1977} written in Greek by Argyros and in a letter Hagler and Haydon sent to Comfort and Negrepontis in 1977. Although Theorem~\ref{double_delta_system} would be the equivalent of letting $\kappa$ equal $\omega$ in \cite[Theorem~1.9, p.~15]{comneg1982} (it should be noted that there are some subtle differences), we decided that it was worth including in this paper a complete proof of this result.

\begin{theorem}\label{double_delta_system} Let $\kappa$ be a cardinal and $\{\kappa_\alpha : \alpha<\cf(\kappa)\}$ be a collection of cardinals with the following properties: 
\begin{enumerate}
\item $\kappa>\cf(\kappa)>\omega$;
\item $\sup\{\kappa_\alpha : \alpha<\cf(\kappa)\}=\kappa$; and
\item for each $\alpha<\cf(\kappa)$, $\kappa_\alpha$ is a regular and uncountable cardinal such that $\max\left\{\cf(\kappa), \sup\{\kappa_\beta : \beta<\alpha\}\right\}<\kappa_\alpha$.
\end{enumerate}

If $\left\{A(\alpha, \gamma) : \alpha<\cf(\kappa)\ \wedge\ \gamma<\kappa_\alpha\right\}$ is a family of finite sets, then there are $m<\omega$, $I\in \left[\cf(\kappa)\right]^{\cf(\kappa)}$, $\{J_\alpha : \alpha\in I\}$, $\{A_\alpha : \alpha\in I\}$ and $A$ which satisfy the following conditions for any distinct $\alpha,\beta\in I$: 
\begin{enumerate}
\item $J_\alpha \in [\kappa_\alpha]^{\kappa_\alpha}$;
\item $A_\alpha \cap A_\beta = A$;
\item if $\gamma,\delta \in J_\alpha$ are distinct, then $A(\alpha, \gamma) \cap A(\alpha, \delta) = A_\alpha$; 
\item whenever $\gamma\in J_\alpha$ and $\delta\in J_\beta$, then $A(\alpha, \gamma) \cap A(\beta, \delta) = A$; and
\item $|A_\alpha|=m$.
\end{enumerate}

\end{theorem}

\begin{proof} Firstly, for each $\alpha<\cf(\kappa)$, Theorem~\ref{delta_system} allows us to find $J_\alpha' \in [\kappa_\alpha]^{\kappa_\alpha}$ and $A_ \alpha$ such that, for any different $\gamma,\delta \in J_\alpha'$, it is satisfied that $A(\alpha,\gamma)\cap A(\alpha,\delta)=A_\alpha$. Next, the inequality $\cf(\kappa)>\omega$ implies that the function $f:\cf(\kappa)\to\omega$ defined by $f(\alpha)=|A_\alpha|$ has a fiber of size $\cf(\kappa)$; in symbols, there exists $m<\omega$ such that $I':=f^{-1}\{m\}$ is an element of $\left[\cf(\kappa)\right]^{\cf(\kappa)}$. Furthermore, another application of Theorem~\ref{delta_system} to the family $\{A_\alpha : \alpha\in I'\}$ allows us to find $I\in [I']^{\cf(\kappa)}$ and $A$ such that $A_\alpha \cap A_\beta = A$, whenever $\alpha, \beta \in I$ are different.

Now let us put $B:=\bigcup\{A_\alpha : \alpha\in I\}$ and recursively define a pair of sequences $\{B_\alpha : \alpha \in I\}$ and $\{J_ \alpha : \alpha\in I\}$ using the following rules: \begin{align*} B_\alpha :&= B\cup \bigcup\left\{A(\beta,\delta) : \beta\in I\ \wedge\ \beta<\alpha\ \wedge\ \delta\in J_\beta\right\} \ \text{and} \\
J_\alpha :&=\left\{\gamma \in J_\alpha' : \left(A(\alpha,\gamma) \setminus A_\alpha\right)\cap B_\alpha =\emptyset\right\}.
\end{align*} Notice that, given $\alpha\in I$, \[|B_\alpha|\leq \left|B\cup \bigcup\{A(\beta,\delta) : \beta<\alpha\ \wedge\ \delta \in J_\beta\}\right| \leq \cf(\kappa)\cdot\sup\{\kappa_\beta : \beta<\alpha\} < \kappa_\alpha.\] For this reason, since $\{A(\alpha,\gamma)\setminus A_\alpha : \gamma\in J_\alpha'\}$ is a disjoint family of cardinality $\kappa_\alpha$, then $|J_\alpha' \setminus J_\alpha|<\kappa_\alpha$ and, therefore, $|J_\alpha|=\kappa_\alpha$.

At this point in the argument it is clear that $m$, $I$, $\{J_\alpha : \alpha\in I\}$, $\{A_\alpha : \alpha\in I\}$ and $A $ satisfy the conditions (1), (2), (3) and (5). To verify (4), let us take distinct $\alpha,\beta\in I$, $\gamma\in J_\alpha$ and $\delta\in J_\beta$, assume without losing generality that $\beta<\alpha$, and concentrate on proving the equality $A(\alpha,\gamma)\cap A(\beta,\delta) = A$. On the one hand, since $\gamma\in J_\alpha$, the relation $\left(A(\alpha,\gamma) \setminus A_\alpha\right)\cap B_\alpha=\emptyset$ implies that $\left( A(\alpha,\gamma) \setminus A_\alpha\right)\cap A(\beta,\delta)=\emptyset$ and thus, $A(\alpha,\gamma)\cap A(\beta ,\delta)\subseteq A_\alpha$. On the other hand, since $\delta\in J_\beta$, the condition $\left(A(\beta,\delta) \setminus A_\beta\right)\cap B_\beta=\emptyset$ tells us that $ \left(A(\beta,\delta) \setminus A_\beta\right)\cap A_\alpha=\emptyset$; therefore $A(\beta,\delta)\cap A_\alpha\subseteq A_\beta$. Hence, the inclusions $A(\alpha,\gamma)\cap A(\beta,\delta )\subseteq A_\alpha$ and $A(\beta,\delta)\cap A_\alpha\subseteq A_\beta$ imply that $A(\alpha,\gamma)\cap A(\beta,\delta)\subseteq A$. Consequently, $A= A_\alpha \cap A_\beta \subseteq A(\alpha,\gamma)\cap A(\beta,\delta) \subseteq A$, i.e., $A(\alpha,\gamma )\cap A(\beta,\delta) = A$.
\end{proof}

\section{Main theorems}\label{secc_main}

For the results that follow we need to make a couple of notational conventions. Let $\{X_\xi : \xi<\lambda\}$ be a family of topological spaces and $A$ be a non-empty subset of $\lambda$. We will use the symbol $X_A$ to denote the topological product $\prod \{X_\xi : \xi \in A\}$. Additionally, $\pi_A : X_\lambda \to X_A$ will be the canonical projection. Lastly, when $A=\{\xi\}$, we will simply write $\pi_\xi$ instead of $\pi_{\{\xi\}}$.

\begin{theorem}\label{central_thm} Let $\{X_\xi : \xi<\lambda\}$ be a collection of topological spaces, $\kappa$ be a cardinal number such that $\kappa>\cf(\kappa)>\omega$ and $\{\kappa_\alpha : \alpha<\cf(\kappa)\}$ be a sequence of cardinals. If for each $\xi<\lambda$, $\kappa$ is a precaliber for $X_\xi$ and $\{\kappa_\alpha : \alpha<\cf(\kappa)\}$ satisfies the $P (\kappa, X_\xi)$ condition, then $\kappa$ is a precaliber for $X_{\lambda}$.

\end{theorem}

\begin{proof} Let $\{U_\gamma : \gamma <\kappa\}$ be a subset of $\tau_{X_\lambda}^+$ formed by canonical open sets. Furthermore, for every $\gamma<\kappa$, let $B_\gamma \in [\lambda]^{<\omega}$ and, for every $\xi\in B_\gamma$, $U_\xi^\gamma \in \tau_{X_\xi}^+$ such that $U_\gamma = \bigcap\{\pi_\xi^{-1}[U_\xi^\gamma] : \xi\in B_\gamma\}$. For the argument to come, we will assume that $\bigcap\{B_\gamma : \gamma <\kappa\} \neq \emptyset$.

Clearly, Theorem~\ref{double_delta_system} allows us to assume, without loss of generality, that there are $m<\omega$, $\{A_\alpha : \alpha<\cf(\kappa)\}$ and $A$ so that the following conditions hold for any distinct $\alpha,\beta<\cf(\kappa)$:
\begin{enumerate}
\item $A_\alpha \cap A_\beta = A$;
\item if $\gamma,\delta \in [\kappa_\alpha,\kappa_{\alpha+1})$ are different, then $B_\gamma \cap B_\delta = A_\alpha$;
\item whenever $\gamma\in [\kappa_\alpha,\kappa_{\alpha+1})$ and $\delta\in [\kappa_\beta,\kappa_{\beta+1})$, then $B_\gamma \cap B_\delta =A$; and
\item $|A_\alpha|=m$.
\end{enumerate}

In particular, items (2) and (3) imply the following property: if $\beta\leq \alpha<\cf(\kappa)$, $\gamma\in [\kappa_\alpha,\kappa_{\alpha +1})$, $\delta\in [\kappa_\beta,\kappa_{\beta+1})$ and $\gamma\neq \delta$, then $\left(B_\gamma\setminus A_\alpha \right) \cap B_\delta = \emptyset$. 

We will divide the rest of the proof into three claims.

\medskip

\noindent {\bf Claim 1.} There exists $J\in [\kappa]^{\kappa}$ such that $\left\{\pi_A [U_\gamma] : \gamma \in J\right\}$ is a centered family.

\medskip

Indeed, since we are assuming that $\bigcap\{B_\gamma : \gamma <\kappa\} \neq \emptyset$, we deduce that $A$ is a non-empty set. Hence, since $\pi_A$ is an open map, the collection $\left\{\pi_A [U_\gamma] : \gamma <\kappa\right\}$ is a subset of $\tau_{X_A}^{+ }$ and thus, Corollary~\ref{cor_2} guarantees the existence of $J\in [\kappa]^{\kappa}$ such that $\left\{\pi_A [U_\gamma] : \gamma \in J \right\}$ is centered.

To simplify the notation we will assume, without losing generality, that $J=\kappa$. In other words, we will suppose from now on the following statement: \begin{align} \left\{\pi_A [U_\gamma] : \gamma <\kappa\right\}\ \text{is a centered family.} \label{statement}
\end{align}

\medskip

\noindent {\bf Claim 2.} There exists a collection $\{J_\alpha : \alpha<\cf(\kappa)\}$ that satisfies the following conditions for each $\alpha<\cf(\kappa)$:

\begin{enumerate}
\item $J_\alpha\subseteq [\kappa_{\alpha+m}, \kappa_{\alpha+m+1})$;  
\item $|J_\alpha|=\kappa_\alpha$; and
\item $\{\pi_{A_{\alpha+m}} [U_\gamma] : \gamma \in J_\alpha\}$ is a centered family.
\end{enumerate}

\medskip

Let $\alpha<\cf(\kappa)$. Given that the collections $\{X_\xi : \xi \in A_{\alpha+m}\}$ and $\{\kappa_{\alpha+l} : 1\leq l\leq m+1\}$ satisfy the conditions of Corollary~\ref{cor_1}, we have that $(\kappa_{\alpha+m+1}, \kappa_{\alpha+1})$ is a precaliber for the space $X_{A_{ \alpha+m}}$. Thus, item (2) of Lemma~\ref{lemma_basic_pairs} ensures that $(\kappa_{\alpha+m+1}, \kappa_\alpha)$ is also a precaliber for $X_{A_{\alpha +m}}$. Therefore, since every element of $\left\{\pi_{A_{\alpha+m}} [U_\gamma] : \gamma \in [\kappa_{\alpha+m}, \kappa_{\alpha+m+1}) \right\}$ is an open and non-empty subset of $X_{A_{\alpha+m}}$, there exists $J_\alpha\in \left[[\kappa_{\alpha+m}, \kappa_ {\alpha+m+1})\right]^{\kappa_\alpha}$ such that $\{\pi_{A_{\alpha+m}} [U_\gamma] : \gamma \in J_\alpha\}$ is a centered family.

Now set $J:=\bigcup\{J_\alpha : \alpha<\cf(\kappa)\}$. Clearly, $|J|=\kappa$. Our last goal is to verify that $\{U_\gamma : \gamma\in J\}$ is a centered family.

Let $I$ be a finite and non-empty subset of $J$. By (\ref{statement}), we can find a point $q \in \bigcap\{\pi_A [U_\gamma] : \gamma \in I\}$. On the other hand, if we let $F:= \{\alpha<\cf(\kappa) : I \cap J_\alpha \neq \emptyset\}$, then for each $\alpha\in F$ we can use claim 2 to fix $r^{\alpha} \in \bigcap\{\pi_{A_{\alpha+m}} [U_\gamma] : \gamma \in I \cap J_\alpha\}$. Next, if for each $\gamma\in J$ we let $\alpha(\gamma)<\cf(\kappa)$ be the only ordinal such that $\gamma \in J_{\alpha(\gamma) }$, then we put \[r:=\left\{\left(\xi,r^{\alpha(\gamma)}(\xi)\right) : \gamma\in I\ \wedge\ \xi\in A_{\alpha(\gamma)+m} \setminus A\right\}.\] Finally, for $\gamma\in I$ and $\xi \in B_\gamma \setminus A_{\alpha(\gamma)+m}$, let us take $p_\xi^\gamma \in U_\xi^{ \gamma}$ and set \[s:= \left\{(\xi,p_\xi^\gamma) : \gamma\in I\ \wedge\ \xi\in B_\gamma \setminus A_{\alpha( \gamma)+m}\right\}.\]

\medskip

\noindent {\bf Claim 3.} $\{q,r,s\}$ is a family of functions with disjoint domains.

\medskip

Evidently, $q$ is a function because it is an element of $X_A$. On the other hand if $\gamma,\delta\in I$ are such that $\alpha(\gamma)\neq \alpha(\delta)$, then $\alpha(\gamma)+m\neq \alpha(\delta)+m$ and so \[ \left( A_{\alpha(\gamma)+m}\setminus A\right) \cap \left(A_{\alpha(\delta)+m}\setminus A\right) = \left( A_{\alpha(\gamma)+m} \cap A_{\alpha(\delta)+m}\right)\setminus A=\emptyset;\] hence $r$ is a function. Lastly, if $\gamma,\delta\in I$ are distinct and such that $\alpha(\gamma)\geq \alpha(\delta)$, then \[\left(B_{\gamma} \setminus A_ {\alpha(\gamma)+m}\right) \cap \left(B_{\delta} \setminus A_{\alpha(\delta)+m}\right) \subseteq \left(B_{\gamma} \setminus A_{\alpha(\gamma)+m}\right) \cap B_{\delta} = \emptyset;\] in particular, $s$ is a function.

Now, from the equalities $\dom(q)=A$, $\dom(r)= \bigcup\left\{A_{\alpha(\gamma)+m} \setminus A : \gamma\in I\right\}$ and $ \dom(s)=\bigcup\left\{B_\gamma\setminus A_{\alpha(\gamma)+m} : \gamma\in I\right\}$ it is easy to see that $\left\{\dom(q),\dom(r),\dom(s)\right\}$ is a disjoint family. Clearly, $\dom(q)$ is disjoint from $\dom(r)\cup\dom(s)$. Additionally, if $\gamma,\delta\in I$ are distinct and we assume, without loss of generality, that $\alpha(\gamma)\geq \alpha(\delta)$, then \[ \left(B_{\gamma } \setminus A_{\alpha(\gamma)+m}\right) \cap \left(A_{\alpha(\delta)+m}\setminus A\right) \subseteq \left(B_{\gamma} \setminus A_{\alpha(\gamma)+m}\right) \cap B_{\delta} = \emptyset.\] Therefore, $\dom(r)\cap \dom(s) =\emptyset$.

Finally, claim 3 implies that $q\cup r \cup s$ can be extended to a point $p\in X_\lambda$. Under these circumstances, it is clear that $p\in\bigcap\{U_\gamma : \gamma \in I\}$. Thus, $\kappa$ is a precaliber for $X_\lambda$.
\end{proof}

\begin{theorem}\label{thm_1} If $X$ is a topological space with precaliber $\kappa$, $\cf(\kappa)>\omega$ and $\lambda$ is a cardinal number, then $X^{\lambda}$ also has precaliber $\kappa$.

\end{theorem}

\begin{proof} When $\kappa$ is regular, Theorem~\ref{thm_Sanin} guarantees that $\kappa$ is a precaliber for $X^{\lambda}$. If $\kappa$ is singular, we use item (2) of Lemma~\ref{lemma_P(kappa,X)} to find a sequence $\{\kappa_\alpha : \alpha<\cf(\kappa)\} $ that satisfies the $P(\kappa,X)$ condition. Finally, Theorem~\ref{central_thm} allows us to conclude that $X^{\lambda}$ has precaliber $\kappa$.
\end{proof}

\begin{theorem}\label{thm_2} Whenever $\{X_\alpha : \alpha<\lambda\}$ is a family of topological spaces, each of them with precaliber $\kappa$, and $\lambda < \cf(\kappa)$, then $X_\lambda$ has precaliber $\kappa$ as well.

\end{theorem}

\begin{proof} For regular $\kappa$, Theorem~\ref{thm_Sanin} tells us that $\kappa$ is a precaliber for $X_\lambda$. If $\kappa$ is singular, our goal is to use the inequality $\lambda<\cf(\kappa)$ to build a sequence like the one needed in Theorem~\ref{central_thm}. Let us fix a strictly increasing sequence of cardinal numbers $\{\lambda_\alpha : \alpha<\cf(\kappa)\}$ such that $\sup\{\lambda_\alpha : \alpha<\cf(\kappa)\}=\kappa$. 

\medskip

\noindent {\bf Claim.} There is a collection $\{\kappa_\alpha : \alpha<\cf(\kappa)\}$ that satisfies the following conditions for each $\alpha<\cf(\kappa)$:

\begin{enumerate}
\item $\max\left\{\cf(\kappa), \lambda_\alpha, \sup\{\kappa_\beta : \beta<\alpha\}\right\} < \kappa_\alpha<\kappa$; 
\item $\kappa_\alpha$ is regular and uncountable; and
\item $(\kappa_{\alpha+1},\kappa_\alpha)$ is a precaliber for $X_\xi$, provided that $\xi<\lambda$.
\end{enumerate}

\medskip

By recursion, suppose that for some $\alpha<\cf(\kappa)$ we have obtained $\{\kappa_\beta : \beta<\alpha\}$ with the desired properties. If we let $\mu:=\max\left\{\cf(\kappa),\lambda_\alpha, \sup\{\kappa_\beta : \beta<\alpha\} \right\}$, we can use item (1) of Lemma~\ref{lemma_P(kappa,X)} to find a cardinal $\mu< \kappa_\alpha^\xi < \kappa$ such that $\left(\kappa_{\alpha}^\xi,\mu\right)$ is a precaliber for $X_\xi$, provided that $\xi<\lambda$. Hence, if we let $\kappa_\alpha := \left(\sup \{\kappa_{\alpha}^\xi : \xi<\lambda\}\right)^{+}$, the relations $\lambda< \cf(\kappa)<\kappa$ imply that $\kappa_\alpha<\kappa$; furthermore, item (2) of Lemma~\ref{lemma_basic_pairs} ensures that, if $\alpha=\gamma+
1$, then $(\kappa_{\gamma+1},\kappa_\gamma)$ is a precaliber for $X_\xi$, whenever $\xi<\lambda$.

Lastly since items (1) to (3) of the claim indicate that $\{\kappa_\alpha : \alpha<\cf(\kappa)\}$ has the properties we need to use Theorem~\ref {central_thm}, we conclude that $X_\lambda$ has precaliber $\kappa$.
\end{proof}

Finally, the combination of Theorems~\ref{thm_Shelah}, \ref{thm_1} and \ref{thm_2} gives us the desired result for singular precalibers and singular calibers:

\begin{theorem}\label{last_thm} The following statements hold for any singular cardinal $\kappa$.

\begin{enumerate}
\item If $X$ is a topological space with precaliber (resp., caliber) $\kappa$, $\cf(\kappa)>\omega$ and $\lambda$ is a cardinal number, then $X^{\lambda}$ also has precaliber (resp., caliber) $\kappa$.
\item Whenever $\{X_\alpha : \alpha<\lambda\}$ is a family of topological spaces, each of them with precaliber (resp., caliber) $\kappa$, and $\lambda < \cf(\kappa)$, then $\prod\{X_\alpha : \alpha<\lambda\}$ has precaliber (resp., caliber) $\kappa$ as well.
\end{enumerate}
\end{theorem}

As a last comment, it is possible to prove a result analogous to Theorem~\ref{last_thm} with respect to the notion of weak precaliber. Recall that a cardinal $\kappa$ is a {\it weak precaliber} for a topological space $X$ if for any set $I$ of cardinality $\kappa$ and for any family $\{U_\gamma : \gamma\in I\} \subseteq \tau_{X}^{+}$, there exists $J\in [I]^{\kappa}$ such that $U_\gamma \cap U_\delta\neq\emptyset$, whenever $\gamma,\delta\in J$. With this concept in mind, each argument in this article can be modified appropriately to prove the following theorem.

\begin{theorem} The following statements hold for any singular cardinal $\kappa$.

\begin{enumerate}
\item If $X$ is a topological space with weak precaliber $\kappa$, $\cf(\kappa)>\omega$ and $\lambda$ is a cardinal number, then $X^{\lambda}$ also has weak precaliber $\kappa$.
\item Whenever $\{X_\alpha : \alpha<\lambda\}$ is a family of topological spaces, each of them with weak precaliber $\kappa$, and $\lambda < \cf(\kappa)$, then $\prod\{X_\alpha : \alpha<\lambda\}$ has weak precaliber $\kappa$ as well.
\end{enumerate}
\end{theorem}

The author wishes to thank the referee for his comments and Dr. \'Angel Tamariz-Mascar\'ua for his invaluable help and direction in the preparation of this article.

\end{document}